\documentclass[13pt]{article}
\usepackage{amsmath}
\usepackage{algorithm}
\usepackage{algorithmic}
\usepackage{amscd}
\usepackage{amsfonts}
\usepackage{amssymb}
\usepackage{amsthm}
\usepackage{epsfig}
\usepackage{slashbox}
\usepackage[all]{xy}
 \newtheorem{lem}{Lemma}[section]

 \newtheorem{prop}[lem]{Proposition}
\newtheorem{prob}[lem]{Problem}
 \theoremstyle{definition}
 
 \theoremstyle{remark}

\begin{document}

\title{AN IMPROVEMENT TO THE  NUMBER FIELD SIEVE}

\author{Qizhi Zhang}

\maketitle

\begin{abstract}
We improve the ``sieve'' part of the  number field sieve used in factoring integer and computing discrete logarithm. The runtime of our method is shorter than that of existing methods. Under some reasonable  assumptions, we prove that it is less than two-thirds of the running time of the algorithm used before  asymptotically with probability greater than 0.6.
\end{abstract}

\section{Introduction}

General number field sieve is used in factoring integer or computing the discrete logarithm. See, for example, \cite{BLP} \cite{Gordon} \cite{Schi1993}.    There are two time consuming parts mainly in the number field sieve. Namely, the part ``sieving'', and the part ``solving the linear equations''. The two parts are relatively independent and have the computational complexity in same order. In \cite{linear}, the authors improved the step ``solving the linear equations'' for discrete logarithm problem. In this paper, we improve the step ``sieve''. Our improvement work for both factoring integer and computing the discrete logarithm. The running time of our algorithm is less than the one in \cite{Schirokauer 2008} \cite{Stevenhagen}  asymptotically. Under some reasonable  assumptions, it is less than $\frac{2}{3}$ of the running time of the algorithm used in  \cite{Schirokauer 2008} \cite{Stevenhagen} asymptotically with probability greater than 0.6.     

In section 2, we give the formulation of the problem which we want to solve, and describe the algorithm used before. In section 3, we describe our algorithm. In section 4, we prove that our algorithm is better than the algorithm used before.

\section{The  problem and conventional algorithm }

Let us consider the following problem: 

\begin{prob}
Let $f(x)$ be a monic polynomial of degree $d$ with integer coefficient that bounded by an integer $m$ and $K$ be an algebraic number field isomorphic to $\mathbb{Q}[x]/(f(x))$. Let $\theta$ be the image of $x$ in $K$ and $Nm : K^\times \longrightarrow \mathbb{Q}^\times$ be the norm map. Let $u$ be a positive integer. We construct a table $T=\{ T(b,a) \}_{0<b \leq u, |a| \leq u }$ of $u$ lines and $2u+1$ cows with
\begin{displaymath}
T(b,a)= \left\{
\begin{array}{ll}
0 & if (a,b)=1 ; \\
(a-bm)Nm(a-b\theta) & if (a,b)\neq 1 .
\end{array}
\right.
\end{displaymath}
Let $y$ be a positive real number called "the smooth bound". For every element in the table, we wish to divide out all of its divisors of the form $l^e$ for all primes $l$ bounded by $y$.
\end{prob}

The most trivial algorithm is the following:

\begin{algorithm}
\caption{Sieve}
\begin{algorithmic}[1]
\FOR{prime $0<l\leq y$, integer $|a|\leq u$, $0<b \leq u$, such that $T(b,a)\neq 0$}
\WHILE{$l|T(b,a)$}
\STATE $T(b,a) \leftarrow T(b,a)/l$
\ENDWHILE
\ENDFOR
\end{algorithmic}
\end{algorithm}

The following improved algorithm is widely used in integer factoring algorithms (see \cite{BLP}, \cite{Stevenhagen})or algorithms of solving the discrete logarithm problem (see  \cite{Gordon} \cite{Schi1993}  \cite{Schirokauer 2008}). 

\begin{algorithm}
\caption{Sieve}
\begin{algorithmic}[1]
\FOR {prime $l \in (0, y]$}
\STATE {$\epsilon_l \leftarrow m \mod l \in \{0, 1, \cdots,  l-1 \}$ }
\STATE {$E_l \leftarrow \{x \in \{0, 1, \cdots, l-1 \}: f(x) \equiv 0 \mod l  \} $}
\ENDFOR
\FOR {integer $b \in (0, u]$}
\FOR {prime $l \in (0, y], l \nmid b$}

\FOR {$a \in [-u,u] \cap (b\epsilon_l +l \mathbb{Z}) $}
\WHILE {$l \mid T(b,a)$}
\STATE {$T(b,a) \leftarrow  T(b,a)/l$}
\ENDWHILE
\ENDFOR

\FOR {$a \in [-u,u] \cap (bE_l +l \mathbb{Z}) $}
\WHILE {$l \mid T(b,a)$}
\STATE {$T(b,a) \leftarrow  T(b,a)/l$}
\ENDWHILE
\ENDFOR

\ENDFOR
\ENDFOR
\end{algorithmic}
\end{algorithm}

In {\bf Algorithm 2}, we do not try to divide all the elements in the table by $l$ more, but divide those divisible by $l$ we know. Then we divide the quotient by $l$ continually as long as it is divisible by $l$.    Roughly speaking, for every $b, l$, we solve the equations   
\begin{displaymath}
a-bm \equiv 0 \mod l \quad \mbox{or       }  Nm(a-b\theta) \equiv 0 \mod l
\end{displaymath} 
of variable $a$, and then sieve.

\section {Our algorithm}
There is unnecessary computing still in algorithm 2. In fact, we can almost know which $T(b, a)$ can be divided by $l$ again, after it divided by $l$ first. Roughly speaking, for every $b, l$, we can almost can solve the equations
\begin{displaymath}
a-bm \equiv 0 \mod l^k \quad \mbox{or       }  Nm(a-b\theta) \equiv 0 \mod l^k
\end{displaymath} 
of variable $a$, for any $k$, and then sieve. Our new algorithm consists of 3 parts. 

{\bf First}, we divide out all $l$-power divisors caused by the term $(a-bm)$: Let $\epsilon _l ^{(1)}$ be the residue for $m$ module $l$. We divide $T(b,a)$ by $l$ and write the quotient in $T(b,a)$ for all $a \in [-u,u] \cap (b \epsilon _l ^{(1)} +l \mathbb{Z})$. Let $\epsilon _l ^{(2)}$ be the residue for $m$ module $l^2$. We divide $T(b,a)$ by $l$ and write the quotient in $T(b,a)$ for all $a \in [-u,u] \cap (b \epsilon _l ^{(2)} +l^2 \mathbb{Z})$. $\cdots$

{\bf Second}, we divide out all the $l$-power divisors caused by the term $Nm (a-b \theta)$, for all $a, b$ such that $b \in (0,u]$ is coprime to $l$, and $a \mod l$ is a single root of the equation $Nm(a- b \theta) \equiv 0 \mod l$. By lemma 3.1 below we can do this as follows:  Let $E _l ^{(1)} \subset \{0, 1, \cdots l-1\}$ be the set of single roots of the equation $f(x)$ module $l$. We can directly compute $E_l ^{(1)}$ by solve equation. We divide $T(b,a)$ by $l$ and write the quotient in $T(b,a)$ for all $a \in [-u,u] \cap (b E _l ^{(1)} +l \mathbb{Z})$. Let $E _l ^{(2)} \subset \{0, 1, \cdots l^2-1 \}$ be the set of single roots of the equation $f(x)$ module $l^2$. We can directly compute $E _l ^{(2)}$ from $E _l ^{(1)}$ by Newton's method. We divide $T(b,a)$ by $l$ and write the quotient in $T(b,a)$ for all $a \in [-u,u] \cap (b E _l ^{(2)} +l^2 \mathbb{Z})$. $\cdots$

{\bf Finally}, we divide out all the $l$-power divisors caused by the term $Nm (a-b \theta)$, for all $a, b$ such that $b \in (0, u]$ is coprime to $l$, and $a \mod l$ is a multiple root of the equation $Nm(a- b \theta) \equiv 0 \mod l$. By lemma 3.1 and lemma 3.2 below we can do this as follows:  Let $\tilde{E} _l ^{(1)} \subset \{0, 1, \cdots l-1 \}$ be the set of multiple roots of the equation $f(x)$ module $l$, we can directly compute $\tilde{E}_l ^{(1)}$ by solving the equation. We divide $T(b,a)$ by $l$ and write the quotient in $T(b,a)$ for all $a \in [-u,u] \cap (b \tilde{E} _l ^{(1)} +l \mathbb{Z})$.  Lemma 3.2 below tells us that whether a root of $f(x)\equiv 0 \mod l$ can be lifted to a root of $f(x) \equiv 0 \mod l^2$ is only dependent on its residue class module $l$. Let $\tilde{E} _l ^{(2)}$ be the subset of $\tilde{E}_l ^{(1)}$ whose elements can be lifted to solutions of $f(x)\equiv 0 \mod l^2$. We can compute $\tilde{E} _l ^{(2)}$ by $\sharp \tilde{E} _l ^{(1)}$'s tests. Then we divide $T(b,a)$ by $l$ and write the quotient in $T(b,a)$ for all $a \in [-u,u] \cap (b \tilde{E} _l ^{(2)} +l\mathbb{Z})$ one after another until $l \nmid T(b,a)$.

Now we give statements and proofs of lemma3.1 and lemma 3.2 mentioned above.

\begin{lem}
If $l \nmid b$, there is a bijective 
\begin{displaymath}
\begin{array}{ccc}
\{x \in \mathbb{Z}/l^e\mathbb{Z}; f(x) \equiv 0 \mod l^e \} & \longrightarrow &  \{ a \in \mathbb{Z}/l^e\mathbb{Z} ; Nm(a-b\theta)\equiv 0 \mod l^e \} \\
x & \mapsto  & bx 
\end{array}
\end{displaymath}
for all $e>0$. Moreover, in the situation $e=1$, the images of simple roots are simple, and the image of multiple roots are multiple. 
\end{lem}

{\bf Proof.} It is because
\[
Nm(a-b \theta)=(-b)^dNm(\frac{a}{b}-\theta)=(-b)^df(\frac{a}{b})
\]
$\hfill \blacksquare$

\begin{lem}
Let $x, y$ be two integers and $f$ be a polynomial over $\mathbb{Z}$. Assume $x \equiv y \mod l$ is a multiple root of $f(x)\equiv 0 \mod l$. Then $x \mod l^2$ is a root of $f(x) \equiv 0 \mod l^2$ if and only if $y \mod l^2$ is a root of $f(x) \equiv 0 \mod l^2$. 
\end{lem}
{\bf Proof.} Let $y=x+kl$ where $k\in\mathbb{Z}$. If $x \mod l^2$ is a root of $f(x)\equiv 0 \mod l^2$, we have
\[
f(y) \in f(x)+f'(x)kl+l^2\mathbb{Z}
\]
by Taylor expansion. On the other hand, we know $f'(x)\equiv 0 \mod l$. Therefore
\[
f(y) \in f(x)+l^2\mathbb{Z}
\] 
$\hfill \blacksquare$

\begin{algorithm}
\caption{sieve}
\begin{algorithmic}[1]
\STATE {{\bf (First)}}
\FOR {prime $l \in (0, y] $}
\FOR {$e =1, 2, \cdots \log _l [u(m+1)]$}
\STATE {$\epsilon _l ^{(e)} \leftarrow m \mod l^e \in \{ 0, 1, \cdots , l^e-1  \} $ }
\ENDFOR
\ENDFOR
\FOR {integer $b \in (0, u]$}
\FOR {prime $l \in (0, y], l \nmid b$}
\FOR {$e =1, 2,  \cdots \log _l [u(m+1)]$}
\FOR {$a \in [-u, u] \cap (b\epsilon _l ^{(e)}+l^e\mathbb{Z}$)}
\STATE {$T(b,a) \rightarrow T(b,a)/l$}
\ENDFOR
\ENDFOR
\ENDFOR
\ENDFOR

\STATE {{\bf (Second)}}
\FOR {prime $l \in (0, y] $}
\FOR {$e=1, 2,  \cdots \log _l [m(d+1)u^d]$}
\STATE $E_l^{(e)} \leftarrow \{ x=0, 1, \cdots l^e-1 : x \mbox{ is a single root of } f(x) \equiv 0 \mod l^e  \}$
\ENDFOR
\ENDFOR
\FOR {integer $b \in (0, u]$}
\FOR {prime $l \in (0, y], l \nmid b$}
\FOR {$e =1, 2, \cdots \log _l [m(d+1)u^d]$}
\FOR {$a \in [-u, u] \cap (bE _l ^{(e)}+l^e\mathbb{Z}$)}
\STATE {$T(b,a) \leftarrow T(b,a)/l$}
\ENDFOR
\ENDFOR 
\ENDFOR
\ENDFOR

\STATE {{\bf (Finally)}}
\FOR {prime $l \in (0, y] $}
\STATE $\tilde{E}_l^{(1)} \leftarrow \{ x=0, 1, \cdots l-1 : x \mbox{ is a multiple root of } f(x) \equiv 0 \mod l  \}$
\STATE $\tilde{E}_l^{(2)} \leftarrow \{ x \in  \tilde{E} _l ^{(1)}: f(x) \equiv 0 \mod l^2 \} $
\ENDFOR

\FOR {integer $b \in (0, u]$}
\FOR {prime $l \in (0, y], l \nmid b$}

\FOR {$a \in [-u, u] \cap (b\tilde{E} _l ^{(1)}+l\mathbb{Z})$}
\STATE $T(b,a) \leftarrow T(b,a)/l$
\ENDFOR

\FOR {$a \in [-u, u] \cap (b\tilde{E} _l ^{(2)}+l\mathbb{Z} )$}
\WHILE {$l \mid T(b,a)$}
\STATE {$T(b,a) \leftarrow T(b,a)/l$}
\ENDWHILE
\ENDFOR 

\ENDFOR
\ENDFOR
\end{algorithmic}
\end{algorithm}

\section {Complexity analysis}
 We will compare the computational complexity of {\bf Algorithm 2} and {\bf Algorithm 3 }. Considering the practical situation, we make the assumption that  $y \leq Ku$ for some constant $K$.  

For $(l, b)=1, e>0$, Let 
\begin{displaymath}
\begin{array}{c}
A_{l}^{b,s} :=\{|a|\leq u; (a,b)=1, a \mbox{ is a single root of }Nm(a-b\theta) \equiv 0 \mod l\} \\
A_{l}^{b,m} :=\{|a|\leq u; (a,b)=1, a \mbox{ is a multiple root of } Nm(a-b\theta) \equiv 0 \mod l\} \\
A_{l^e} ^{b,s} := \{a \in A_{l}^{b,s}; Nm(a-b\theta) \equiv 0 \mod l^e  \}   \\
A_{l^e} ^{b,m}:= \{a \in A_l ^{b,m}; Nm(a-b\theta) \equiv 0 \mod l^e \}\\
A_{l^e}^b:=\{|a|\leq u; (a,b)=1,  Nm(a-b\theta) \equiv 0 \mod l^e\}  \\
B _{l^e} ^b:= \{|a|\leq u; (a,b)=1, a-bm \equiv 0 \mod l^e \}  
\end{array}
\end{displaymath}

In Algorithm 2, the complexity of line $1$ --line $4$ is an infinitesimal of the complexity of line $5$--line$18$ as $u \rightarrow \infty$. From line $5$, the complexity of sieving the elements in the $b$-th line of the table by prime $l$ is
\begin{displaymath}
\begin{array}{rll}
C_l^b=& \sharp B_l^b & (T(b,a) \leftarrow T(b,a)/l \mbox{ for } a \in B_l ^b ) \\
&+ \sharp B_l^b & (\mbox{ try to divide } T(b,a) \mbox{ by } l \mbox{ for } a \in B_l ^b \mbox{ again })\\
&+ \sharp B_{l^2}^b & (\mbox{ try to divide } T(b,a) \mbox{ by } l \mbox{ for } a \in B_{l^2}^b ) \\
&+ \sharp B_{l^3}^b & (\mbox{ try to divide } T(b,a) \mbox{ by } l \mbox{ for } a \in B_{l^3}^b) \\
&+ \cdots & \\
&+ \sharp A_l^b & (\mbox{ try to divide } T(b,a) \mbox{ by } l \mbox{ for } a \in A_l ^b) \\ 
&+ \sharp A_l^b & (\mbox{ try to divide } T(b,a) \mbox{ by } l \mbox{ for a } \in A_l^b \mbox{ again})\\
&+ \sharp A_{l^2}^b & (\mbox{ try to divide } T(b,a) \mbox{ by } l \mbox{ for } a \in A_{l^2} ^b)\\
&+ \sharp A_{l^3}^b & (\mbox{ try to divide } T(b,a) \mbox{ by } l \mbox{ for } a \in A_{l^3}^b)\\
&+ \cdots & \\  
=& \sharp B_l^b + \sharp B_l^b + \sharp B_{l^2}^b + \sharp B_{l^3}^b + \cdots &\\
& + \sharp A_l^{b,s} +  \sharp A_l^{b,s} + \sharp A_{l^2}^{b,s} +\sharp A_{l^3}^{b,s} + \cdots& \\
& + \sharp A_l^{b,m} + \sharp A_l^{b,m} + \sharp A_{l^2}^{b,m} +\sharp A_{l^3}^{b,m} + \cdots &\\
=& \sharp B_l^b (1+1+\frac{1}{l}+\frac{1}{l^2}+ \cdots) & \\
&+ \sharp A_l ^{b,s}(1+1+\frac{1}{l}+\frac{1}{l^2}+ \cdots) & \\
& + \sharp A_l^{b,m} + \sharp A_l^{b,m} + \sharp A_{l^2}^{b,m} +\sharp A_{l^3}^{b,m} + \cdots &\\
=& \frac{2l-1}{l-1} \sharp B_l^b + \frac{2l-1}{l-1} \sharp A_l^{b,s} & \\
& + 2 \sharp A_l^{b,m} + \sharp A_{l^2}^{b,m} +\sharp A_{l^3}^{b,m} + \cdots &\\

\end{array}
\end{displaymath}  
Therefore the total complexity of Algorithm 2 is 
\begin{displaymath}
\begin{array}{rl}
(1+o(1))\sum _{\mbox{integer } b \in [1, u]} \sum _{\mbox{prime } l \in [2, y], l\nmid b} C_l^b \quad \mbox{ as } u \rightarrow \infty
\end{array}
\end{displaymath}

In the first part of Algorithm 3, the complexity of line $2$--line $6$ is an infinitesimal of the complexity of line $7$--line $15$ as $u \rightarrow \infty$. The latter is
\begin{displaymath}
\begin{array}{rl}
& \sharp B_l^b +  \sharp B_{l^2}^b + \sharp B_{l^3}^b + \cdots \\
= & \sharp B_l^b (1+\frac{1}{l}+\frac{1}{l^2})  \\
=&\frac{l}{l-1}\sharp B_l^b \\
\end{array}
\end{displaymath}
 The complexity of line $17$--line $21$ is an infinitesimal of the complexity of line $22$--line $30$ as $u\rightarrow \infty$. The latter is
\begin{displaymath}
\begin{array}{rl}
& \sharp A_l^{b,s} +  \sharp A_{l^2}^{b,s} + \sharp A_{l^3}^{b,s} + \cdots \\
= & \sharp A_l^{b,s} (1+\frac{1}{l}+\frac{1}{l^2})  \\
=&\frac{l}{l-1}\sharp A_l^{b,s} \\
\end{array}
\end{displaymath}
The complexity of line $32$--line $35$ is an infinitesimal of the complexity of line $36$--line $47$ as $u \rightarrow \infty$. The latter is
\begin{displaymath}
 \sharp A_l^{b,m} +  \sharp A_{l^2}^{b,m} + \sharp A_{l^3}^{b,m} + \cdots 
\end{displaymath}

Therefore the total complexity of Algorithm 3 is 
\begin{displaymath}
(1+o(1))\sum _{\mbox{integer } b \in [1, u]} \sum _{\mbox{prime } l \in [2, y], l\nmid b} D_l^b \quad \mbox{ as }u\rightarrow \infty ,
\end{displaymath}
where
\begin{displaymath}
D_l^b=\frac{l}{l-1}\sharp B_l^b+\frac{l}{l-1} \sharp A_l^{b,s}+\sharp A_l^{b,m} +  \sharp A_{l^2}^{b,m} + \sharp A_{l^3}^{b,m} + \cdots 
\end{displaymath}

It is easy to see that the complexity of Algorithm 3 is less than the complexity of Algorithm 2 gradually. Moreover, if for any $(l, b)=1$, we have $\sharp A_{l^2}^{b,m}=0$ and then we have $D_l^b <\frac{2}{3}C_l^b$, because
\begin{displaymath}
\begin{array}{l}
\frac{l}{l-1}\sharp B_l^b  \leq \frac{2}{3} \times \frac{2l-1}{l-1} \sharp B_l^b \quad \mbox{ for any } (l, b)=1 \\
\frac{l}{l-1} \sharp A_l^{b,s} \leq \frac{2}{3} \times \frac{2l-1}{l-1} \sharp A_l^{b,s} \quad \mbox{ for any } (l, b)=1\\
\sharp A_l^{b,m} < \frac{2}{3} \times 2\sharp A_l^{b,m}  \quad \mbox{ for any } (l, b)=1\\
\sharp A_{l^e} ^{b,m}=0 \quad \mbox{ for any } (l, b)=1 \mbox{ and all } e>1
\end{array}
\end{displaymath} 
Therefore we get
\begin{prop}
Let $K>0$ be a constant. Let $u \rightarrow \infty$ and $y<Ku$. Then the complexity of Algorithm 3 is less than the complexity of Algorithm 2 asymptotically. Moreover, if for any $(l, b)=1$ we have $\sharp A_{l^2}^{b,m}=0$, then the complexity of Algorithm 3 is less than  $\frac{2}{3}$ of the complexity of Algorithm 2  asymptotically.          $\hfill \blacksquare$
\end{prop}

The following proposition tolls us that the condition "for any $(b, l)=1, \sharp A_{l^2}^{b,m}=0$" has much chance to be realized.

\begin{prop}
 Suppose $f(x)$ is a random polynomial of degree $d$ over $\mathbb{Z}$ such that $f(x) \mod l^2$ is uniform distribution on  $\{ h(x) \in  \mathbb{Z}/l^2\mathbb{Z}[x]; \deg h \leq d \}$ for all prime $l \leq y$, and $\{R_l=\phi \}_{\mbox{prime } l \leq y}$ are independent random events, where $R_l:=\{x \in \mathbb{Z}/l^2\mathbb{Z}; f(x)\equiv 0 \mod l^2, x \mbox{ is a multiple root of } f(x)\equiv 0 \mod l \}$ for any prime $l \leq y$. Then the probability of event $(\sharp A_{l^2}^{b,m}=0, \mbox{ for any }b\in (0, u], \mbox{ prime }l \leq y, \mbox{ s.t } (l, b)=1)$ is greater than $0.6$.
\end{prop} 

{\bf Proof.} For any prime $l$, we have
\begin{displaymath}
\begin{array}{rll}
&P(\sharp R_l \neq 0) & \\
=& P(\exists \quad  i=0, 1, \cdots l-1, \quad  s.t. \quad  i \in R_l) & \mbox{(lemma 3.2)} \\
\leq & \sum _{i=0} ^{l-1} P(i \in R_l) &   
\end{array}
\end{displaymath}
For any $i=0, 1, \cdots l-1$, the asuumption that  $f(x) \mod l^2$ is uniform distribution on  $\{ h(x) \in  \mathbb{Z}/l^2\mathbb{Z}[x]; \deg h \leq d \}$ and lemma 4.3 below show

\[
P(f(i)\equiv 0 \mod l^2, \quad f'(i)\equiv 0 \mod l)=\frac{1}{l^3}
\]
i.e. $P(i \in R_l)=\frac{1}{l^3}$.  Hence $P(\sharp R_l \neq 0) \leq \frac{1}{l^2}$. The assumption that  $\{R_l=\phi \}_{\mbox{prime } l \leq y}$ are independent random events implies
\begin{displaymath}
\begin{array}{rl}
 & P(\sum _{\mbox{prime } l \leq y} \sharp R_l=0) \\
=& \prod _{\mbox{prime } l \leq y } P(\sharp R_l =0) \\
=& \prod _{\mbox{prime } l \leq y } (1-\frac{1}{l^2})\\
>& \prod _{l: \mbox{ prime }}(1-\frac{1}{l^2}) \\
=&1/\zeta(2)>0.6 ,
\end{array}
\end{displaymath}
where $\zeta(s)$ is the Riemann's Zeta function.                              

Let $R_l^b:=\{a \in \mathbb{Z}/l^2\mathbb{Z} ; Nm(a-b \theta) \equiv 0 \mod l^2, a \mbox{ is a multiple root of } Nm(a-b\theta) \equiv 0 \mod l \}$. From lemma 3.1, we know 
\begin{displaymath}
\sharp R_l^b =\sharp R_l \quad \mbox{ for all } (b, l)=1
\end{displaymath}
Therefore
\begin{displaymath}
\begin{array}{rl}
&P(\sharp A_{l^2}^{b,m}=0, \mbox{ for any }b\in (0, u], \mbox{ prime }l \leq y, \mbox{ s.t } (l, b)=1)\\
\geq &P(\sharp R_l^b=0, \mbox{ for any }b\in (0, u], \mbox{ prime }l \leq y, \mbox{ s.t } (l, b)=1)\\
=& P(\sharp R_l=0, \mbox{ for any  prime }l \leq y)>0.6
\end{array} 
\end{displaymath}   $\hfill \blacksquare$

Now we give the statement and proof of lemma 4.3 mentioned above.

\begin{lem}
Let $l$ be a prime and $d$ be a positive integer. For any $i\in \mathbb{Z}/l^2\mathbb{Z}$, we have
\begin{displaymath}
P( h(i)\equiv 0 \mod l^2, h'(i)\equiv 0 \mod l \mid h(x) \in \mathbb{Z}/l^2\mathbb{Z}[x], \deg h=d, \mbox{monic})=\frac{1}{l^3}
\end{displaymath}
\end{lem}
{\bf Proof.}Consider the surjective homomorphism of abelian group

\begin{displaymath}
\begin{array}{ccl}
\{h(x)\in \mathbb{Z}/l^2\mathbb{Z}[x]; \deg h \leq d\}& \longrightarrow &\mathbb{Z}/l^2\mathbb{Z} \oplus \mathbb{Z}/l\mathbb{Z} \\
h(x)& \mapsto & (h(i) \mod l^2, h'(i) \mod l).
\end{array}
\end{displaymath}
We have
\begin{displaymath}
P( h(i)\equiv 0 \mod l^2, h'(i)\equiv 0 \mod l \mid h(x)\in \mathbb{Z}/l^2\mathbb{Z}[x], \deg h \leq d ) = \frac{1}{l^3} .
\end{displaymath}
Similarly, we have
\begin{displaymath}
P( h(i)\equiv 0 \mod l^2, h'(i)\equiv 0 \mod l \mid h(x)\in \mathbb{Z}/l^2\mathbb{Z}[x], \deg h \leq d-1 ) = \frac{1}{l^3} .
\end{displaymath}

Let
\begin{displaymath}
\begin{array}{l}
H:=\{ h \in \mathbb{Z}/l^2\mathbb{Z}[x] ; \deg f \leq d \}  \\
H^0:=\{ h \in H ; h(i)\equiv 0 \mod l^2, h'(i) \equiv 0 \mod l \} \\
H_{l\mathbb{Z}}:=\{ h \in H ; \mbox{ the leading coefficient of } h \mbox{ is in }  l\mathbb{Z} \}   \\
H_{l\mathbb{Z}}^0:=H^0 \cap H_{l\mathbb{Z}} \\
H_c:=\{ h \in H ; \mbox{ the leading coefficient of } h \mbox{ is } c \} \quad \mbox{ for any } c \in \mathbb{Z}/l^2\mathbb{Z}\\
H_c^0:=H_c \cap H^0  \quad \mbox{ for any } c \in \mathbb{Z}/l^2\mathbb{Z}  
\end{array}
\end{displaymath}
It is easy to see that
\begin{displaymath}
\begin{array}{l}
H_0=\{ h \in \mathbb{Z}/l^2\mathbb{Z}[x] ; \deg g \leq d-1 \} \\
H_0^0=\{ h \in H_0 ; h(i) \equiv 0 \mod l^2, h'(i) \equiv 0 \mod l \}
\end{array}
\end{displaymath}
Hence, we have
\begin{displaymath}
P( h\in H_0 ^0 \mid h \in H_0)=P(h \in H^0 \mid h\in H)=\frac{1}{l^3}
\end{displaymath}

Let us consider the commutative diagram of abelian groups 
\begin{displaymath}
\xymatrix@!=1.5pc{
                  0 \ar[rr] && H_{l\mathbb{Z}}^0 \ar[rr] \ar@{_{(}->}[d] &&
                  H^0 \ar@{_{(}->}[d] \ar[rr] && \mathbb{Z}/l^2\mathbb{Z} \ar[rr] \ar[d]^{id} && 0 \\
                  0 \ar[rr] && H_{l\mathbb{Z}}  \ar[rr]  && H \ar[rr]  && \mathbb{Z}/l^2\mathbb{Z} \ar[rr] && 0 \\
}
\end{displaymath}
where the map from $H$ to $\mathbb{Z}/l^2\mathbb{Z}$ is defined by $h \mapsto  (\mbox{ the leading coefficient of }h)$. The vertical map in the right side is an identity, hence we have
\[
P( h \in H_{l\mathbb{Z}}^0 \mid h \in H_{l\mathbb{Z}}) = P ( h \in H^0 \mid h \in H  )= \frac{1}{l^3}
\]  
On the other hand,
\begin{displaymath}
\begin{array}{rl}
&P( h\in H_{l\mathbb{Z}}^0 \mid h\in H_{l\mathbb{Z}} )  \\
=& \sum _{c=0} ^{l-1} P( h\in H_{l\mathbb{Z}}^0 \mid h \in H_{cl} ) P( h \in H_{cl} \mid h\in H_{l\mathbb{Z}}  ) \\
=& \sum _{c=0} ^{l-1} P( h\in H_{cl}^0 \mid h \in H_{cl} ) P( h \in H_{cl} \mid h\in H_{l\mathbb{Z}}  ) \\
=& [P ( h \in H_0^0 \mid h\in H_0)+ \sum _{c=1}^{l-1}] \times \frac{1}{l}  \\
=&\frac{1}{l} [\frac{1}{l^3}+ \sum _{c=1} ^{l-1}P( h \in H_{cl}^0 \mid h \in H_{cl}) ]
\end{array}
\end{displaymath}
Hence we have
\[
\sum_{c=1}^{l-1} P (  h\in H_{cl}^0 \mid h \in H_{cl}) = \frac{l-1}{l^3} 
\]
For any $c \in (\mathbb{Z}/l^2\mathbb{Z})^\times$, we have a commutative diagram of sets
\begin{displaymath}
\xymatrix@!=1.5pc{
                   H_1^0 \ar[rr]^\sim  \ar@{_{(}->}[d] &&  H^0_c  \ar@{_{(}->}[d]  \\
                   H_1  \ar[rr]^\sim  && H_c
}
\end{displaymath}
where the horizontal map is defined by $h \mapsto ch$. Hence we have
\[
P( h\in H_1^0 \mid h \in H_1  )=P( h \in H_c^0 \mid h \in H_c ) \quad \mbox{ for any }c \in (\mathbb{Z}/l^2\mathbb{Z})^\times
\]
Therefore
\begin{displaymath}
\begin{array}{rl}
\frac{1}{l^3}=&P( h\in H^0 \mid h\in H)\\
=& \sum _{c \in \mathbb{Z}/l^2\mathbb{Z}}P( h \in H^0 \mid h \in H_c )P( h\in H_c \mid H) \\
=& \sum _{c \in \mathbb{Z}/l^2\mathbb{Z}}P( h\in H_c^0 \mid h\in H_c)P(H_c \mid H) \\
=&[P( h \in H_0^0 \mid h \in H_0)+ \sum _{c=1}^{l-1}P( h\in H_{cl}^0 \mid h\in H_{cl})+ \sum _{c \in (\mathbb{Z}/l^2\mathbb{Z})^\times}P( h \in H_c^0 \mid h\in H_c) ]\times \frac{1}{l^2} \\
=&[ \frac{1}{l^3}+ \frac{l-1}{l^3}+ (l^2-l)P( h \in H_1^0 \mid h\in H_1) ]\times \frac{1}{l^2}
\end{array}
\end{displaymath}
Hence we have
\[
P( h \in H_1 ^0 \mid h\in H_1)=\frac{1}{l^3}            
\]
$\hfill \blacksquare$

Finally, from proposition $4.1$ and proposition $4.2$, we get the main conclusion of this paper:
\begin{prop}
Let $K>0$ be a constant. Let $u \rightarrow \infty$ and $y<Ku$, then the complexity of Algorithm 3 is less than the complexity of Algorithm 2 asymptotically. Moreover,  suppose $f(x)$ is a random polynomial of degree $d$ over $\mathbb{Z}$ such that $f(x) \mod l^2$ is uniform distribution on  $\{ h(x) \in  \mathbb{Z}/l^2\mathbb{Z}[x]; \deg h \leq d \}$ for all prime $l \leq y$, and $\{R_l=\phi \}_{\mbox{prime } l \leq y}$ are independent random events, where $R_l:=\{x \in \mathbb{Z}/l^2\mathbb{Z}; f(x)\equiv 0 \mod l^2, x \mbox{ is a multiple root of } f(x)\equiv 0 \mod l \}$ for any prime $l \leq y$, then the complexity of Algorithm 3 is less than $\frac{2}{3}$ of the complexity of Algorithm 2  asymptotically with probability greater than 0.6.          $\hfill \blacksquare$
\end{prop}

\noindent $\mathbf{Acknowledgement.}$ I would like to thank professor Takeshi
Saito, who gave me much valuable advice.

\end{document}